\documentclass[12pt]{article}
\usepackage{fullpage,amsfonts,amsmath,amsthm}

\def\QQ{\mathbb{Q}}
\def\ZZ{\mathbb{Z}}

\newtheorem*{theorem*}{Theorem}
\newtheorem*{prop*}{Proposition}

\begin{document}

\title{In search of Robbins stability}
\author{Kiran S. Kedlaya \\
Department of Mathematics \\ Massachusetts Institute of Technology 
\\ 77 Massachusetts Avenue
\\ Cambridge, MA 02139 \\ \texttt{kedlaya@math.mit.edu} \\ \\
James G. Propp \\
Department of Mathematics \\ 
University of Wisconsin \\ Madison, WI 53706 \\ \texttt{propp@math.wisc.edu}}
\date{September 27, 2004}

\maketitle

\begin{abstract}
We speculate on whether a certain $p$-adic stability 
phenomenon, observed by David Robbins empirically for Dodgson 
condensation, appears in other nonlinear recurrence relations that 
``unexpectedly'' produce integer or nearly-integer sequences. We exhibit 
an example (number friezes) where this phenomenon provably occurs.
\end{abstract}

This note may be viewed as an addendum to Robbins's note \cite{robbins}
in this volume.
Its purpose is to speculate on whether the $p$-adic stability phenomenon
that Robbins observed empirically for Dodgson condensation appears in other
nonlinear recurrence relations that ``unexpectedly'' produce integer or 
nearly-integer
sequences, and to provide an example where this provably occurs.

In order to carry out this speculation, we'll phrase Robbins's observation
in a somewhat more general framework. For us, a \emph{recurrence relation}
over a field $K$ will consist of a finite partially ordered set $S$ plus,
for each $s \in S$, a rational function $f_s = P_s/Q_s$ 
over $K$ in the indeterminate vector $(x_t)_{t<s}$.
(The restriction to $S$ finite does not concede any generality for
our purposes: to consider an infinite
recurrence, look instead at all of its finite truncations.) 
We also assume (for simplicity)
that the partial order on $S$ is generated by the relation
in which $t$ is less than 
$s$ if $f_s$ is nonconstant as a function of $x_t$ alone.
In this case, $s \in S$ is minimal for the partial order if and only if 
$f_s$ is a constant function; we thus use the term \emph{initial}
interchangeably with ``minimal''.

Before proceeding further, it will be helpful to set up some more notation.
For $I = (i_s)_{s \in S}$ 
a tuple of nonnegative integers,
we write $x^I$ for $\prod_{s \in S} x_s^{i_s}$; for any
function $g: S \to K$, we write $g^I$ for $\prod_{s \in S} g(s)^{i_s}$.
Write $P_s = \sum_I a_{s,I} x^I$ and $Q_s = \sum_I b_{s,I} x^I$, where
$a_{s,I}$ and $b_{s,I}$ are zero for all but finitely many $I$,
and $P_s$ and $Q_s$ have no common polynomial factor.

Suppose now that $K$ is equipped with a discrete (nonarchimedean) valuation
$v$, e.g., $K = \QQ$ with the $p$-adic valuation for some prime $p$.
Suppose also that the $P_s$ and $Q_s$ are normalized so that
$v(a_{s,I}) \geq 0$ and $v(b_{s,I}) \geq 0$ for all $s$ and $I$,
and so that for each $s$,
\[
\min_I \{\min\{v(a_{s,I}), v(b_{s,I})\}\} = 0.
\]
Suppose further that there exists a function
$g: S \to K$ such that $g(s) = f_s(g)$ for all
$s \in S$; note that $g$ is unique if it exists, and the only obstruction
to its existence is the vanishing of $Q_s$ for some $s$.
That is, $g$ is the unique solution of the recurrence, and satisfies
\[
g(s) = \frac{\sum_I a_{s,I} g^I}{\sum_I b_{s,I} g^I}
\]
for all $s \in S$.

Now fix a positive integer $N$. We denote by $*$ any element of $K$
with $v(*) \geq N$;
here we intend that two different occurrences of $*$ may refer to two
different numbers. With this convention, we have the following
simplification rules:
\begin{align*}
* + * &= * \\
(1 + *)(1 + *) &= 1 + * \\
(1 + *)/(1 + *) &= 1 + *.
\end{align*}
We also have $c * = *$ whenever $v(c) \geq 0$.

Define an \emph{$N$-perturbation} of the recurrence as any function $g': 
S \to K$ such that for each $s \in S$,
\[
g'(s) = \frac{\sum_I (1+*) a_{s,I} (g')^I}{\sum_I (1+*)
b_{s,I} (g')^I}.
\]
In case $s$ is initial, this yields $g'(s) = g(s)(1 + *)$;
this is the same as saying that $v(g'(s) - g(s)) \geq v(g(s)) + N$.

The point of this definition is that, in the case $K = \QQ_p$, 
$g'$ is a possible result of
computing $f_s(g')$ using $p$-adic floating point numbers with $N$-digit
mantissas. Specifically, recall from \cite{robbins} that a ``$p$-adic
floating point number with an $N$-digit mantissa''
consists of a pair $(a,e)$, where
the ``mantissa'' $a$ is an invertible element of $\ZZ/p^n\ZZ$ and
the ``exponent'' $e$ is any integer.
This pair is used to represent any $p$-adic number $\tilde{a} p^e$ such that
$\tilde{a}$ is invertible in $\ZZ_p$ and the image of 
$\tilde{a}$ under the natural map from $\ZZ_p$ to 
$\ZZ/p^n \ZZ$ is $a$. Hence two numbers $r$ and $s$ admit
the same representation if and only if $r = s(1 + p^N u)$ for some
$u \in \ZZ_p$, i.e., if $v(s/r - 1) \geq N$.

One can then reimagine $p$-adic floating point arithmetic as being
carried out with actual $p$-adic numbers, except that at any point
in an arithmetic operation, a gremlin may come along and multiply any
value by a factor of the form $1 + *$. In this interpretation,
$g'(s)$ is then allowed to be any result of computing $f_s(g')$
in the presence of such gremlins. (Note that any ``gremlin factor''
applied after adding two numbers together can be absorbed into the
gremlin factors by which each summand is multiplied. Also, the
reciprocal of a gremlin factor is itself a gremlin factor.)

Given an $N$-perturbation $g'$, define its
\emph{projected precision loss} $r_s(g')$
at $s \in S$ as
\[
r_s(g') = \max_{t \leq s} \{v(Q_t(g'))\};
\]
this generalizes the notion of ``condensation error'' introduced by Robbins.
Note that the projected precision loss is determined by the
\emph{computed} denominators rather than the \emph{actual} denominators,
which would be the $v(Q_t(g))$; these often but do not always coincide.
Note also that $r_s(g') = 0$ when $s$ is initial (because
the only term in the maximum is $v(Q_s(g')) = v(1) = 0$),
and that $r_s(g') \geq r_t(g')$ whenever $t \leq s$, i.e., the
bound gets larger (i.e., worse) as you go along.

We say that the recurrence  \emph{exhibits Robbins stability} if for any 
positive integer $N$, any 
$N$-perturbation $g'$, and any $s \in S$, if
$r_s(g') < N$, then
\[
v(g'(s) - g(s)) \geq N - r_s(g') + \min\{0, v(g(s))\}.
\]
Robbins's conjecture in \cite{robbins}, made on the basis of
copious numerical evidence, then essentially (but see next
paragraph) amounts to the statement that
the recurrence obtained from Dodgson condensation of a matrix of 
indeterminates (indexed by the connected minors) exhibits
Robbins stability. (Note that 
the term $\min\{0, v(g(s))\}$ drops out in Robbins's
case because $v(g(s))$ is always nonnegative; this seems to be
warranted by experimental evidence, as we note at the very end.)

It may be more accurate to speak here of ``weak Robbins stability'',
as we are actually generalizing a slightly restricted version
of Robbins's conjecture. That is because Robbins permits the ``borderline''
case $r_s(g') = N$; indeed, the comment in \cite{robbins}
that ``a quantity can be accurate to zero places'' suggests that 
this permission was deliberate. However, we are not entirely sure
whether to believe the conjecture in the borderline case, and
our proof of Robbins stability
in the one nontrivial case we can prove (see the Theorem
below) does not handle the borderline case; a future clarification
of this issue would be welcome.

It may be helpful to recall (or rather, to attempt to reconstruct)
some of Robbins's motivation for making his original conjecture.
The quantity $N - r_s(g')$ measures the extent to which we can
distinguish the denominators we have divided by so far from zero.
To the extent that we can make this distinction, we expect that
Dodgson condensation should continue to work; this expectation
is formalized in the inequality defining the stability property.

However, the assertion that $N - r_s(g')$ measures our ability
to distinguish denominators from zeroes is only really meaningful
if those denominators are integral. This suggests that it may not
be wise to expect stability for recurrences in which denominators
occur in an unsystematic fashion; this caution is borne out
by a simple example, which we now give.

Take $S = \{0,1,2,3,4,5,6,7\}$, equipped with the ordering that agrees
with the usual ordering except that $0$ and $1$ are not comparable, 
and consider
the recurrence over $\QQ$ given by
\[
x_0 = 5, \quad
 x_1 = -5, \quad x_{n} = \frac{x_{n-1} - 1}{x_{n-2}} \qquad (n=2,\dots,7).
\]
The function $g$ in this case takes the values
\[
5, -5, -\frac{6}{5}, \frac{11}{25}, \frac{7}{15}, -\frac{40}{33},
-\frac{365}{77}, \frac{663}{140}.
\]
Let $v$ denote the $2$-adic valuation; then the function
$g'$ taking the values
\[
5, -5, -\frac{6}{5}, \frac{11}{25}, -\frac{793}{15}, -\frac{4040}{33},
\frac{20365}{8723}, -\frac{17463}{1601860}
\]
is an $N$-perturbation for $N=6$, because
\[
g'(4) = \frac{11/25 - (1 - 2^6)}{-6/5}
\]
and $g'(n) = f_n(g')$ for $n=5,6,7$. 
The projected precision loss is
\[
r_7(g') = \max\{v(5), v(-5), v(-6/5), v(11/25), v(7/15), v(-40/33)\}
= 3,
\]
and $v(663/140) = -2$, so Robbins stability would predict that
\[
v(-17463/1601860 - 663/140) \geq N - r_7(g') + \min\{0, v(663/140)\}
= 6 - 3 - 2 = 1.
\]
However, $-17463/1601860 - 663/140 = -2661195/560651$ has valuation 0,
so the recurrence does not exhibit Robbins stability.

As noted before, it is unclear whether one should expect Robbins
stability to be exhibited by recurrences with ``unpredictable''
denominators. However, there is a wide class of recurrences in 
which denominators either do not occur, or occur in a limited
and systematic fashion; these are the recurrences which exhibit
the ``Laurent phenomenon'', in the parlance of Fomin and Zelevinsky
\cite{fz}.
That paper establishes that a number of interesting
recurrences (like Dodgson condensation)
have the following property: if one views the initial constants
as distinct indeterminates, the noninitial terms turn out to be polynomials
in these indeterminates and their inverses.
(See \cite{propp} for an online discussion of such recurrences
and related topics.)

Among recurrences admitting the Laurent phenomenon,
Dodgson condensation is but one example, and it seems (to us, anyway)
that the unexpected cancellations that contribute to the Laurent 
phenomenon may in the condensation case must have something to do with
the unexpectedly strong bound on the precision loss predicted by
Robbins stability. We thus pose the question: do other Laurent
recurrences exhibit Robbins stability?

One can trivially construct many recurrences exhibiting Robbins
stability, by considering those
for which $Q_s = 1$ for all $s$,
so that no divisions are ever performed in the calculation
and hence $r_s(g') = 0$ for all $s \in S$.
In fact, these recurrences have a much stronger property.
\begin{prop*}
Suppose $Q_s = 1$ for all $s$. Then
for any $N$-perturbation $g'$,
$v(g'(s) - g(s)) \geq N$ (and hence $v(g'(s)) \geq 0$)
for all $s \in S$.
\end{prop*}
\begin{proof}
We proceed by induction on $s$; for $s$
minimal, the desired inequality is given directly by the definition
of an $N$-perturbation, so we assume that $s$ is nonminimal and that
\[
g'(t) = g(t) + * \qquad \mbox{for all $t < s$.}
\]
In particular, $v(g'(t)) \geq 0$ for all $t <s$.

We now begin a second induction to show that
$(g')^I = g^I + *$ for all tuples $I$ of nonnegative integers
indexed by the set of $t \in S$ with $t<s$;
this induction will be on the sum of
the entries of $I$. If this sum is zero, then the desired equality
is the trivially true $1 = 1 + *$. Otherwise, given a tuple $I$
for which the claim is known for all tuples of smaller sum,
choose some $t$ at which $I$ has a nonzero component, and let $J$
be the tuple obtained by decreasing this component by 1. Then
$g^I = g^J g(t)$ and likewise for $g'$, $(g')^{J} = g^J + *$
by the inner induction hypothesis, and $g'(t) = g(t) + *$ by the
outer induction hypothesis. These imply that $g'(t)$ and $(g')^{J}$ have
nonnegative valuation, and so
\begin{align*}
(g')^I &= (g')^J g'(t) \\
&= (g^J + *)(g(t) + *) \\
&= g^J g(t) + g(t) * + g^J * + \, * \\
&= g^J g(t) + * \\
&= g^I + *.
\end{align*}
This completes the inner induction, so we may conclude that
$(g')^I = g^I + *$ for all $I$.

To complete the outer induction, note that
\begin{align*}
g'(s) - g(s) &= \sum_I (a_{s,I} + *) (g')^I - a_{s,I} g^I \\
&= \sum_I (g')^I * - \sum_I a_{s,I} ((g')^I - g^I) \\
&= \sum_I * - \sum_I a_{s,I} * = *
\end{align*}
since $v(a_{s,I}) \geq 0$ by hypothesis.
\end{proof}

On the other hand, it seems not so easy to establish  that
Robbins stability is exhibited by any recurrences, even ones
exhibiting the Laurent phenomenon, in which nontrivial divisions
take place. However, we have succeeded in doing so in one case,
which we now describe; it is a form of a recurrence 
of Conway and Coxeter \cite{conway}, which we will refer to here
as the ``number frieze'' recurrence.

Fix a positive integer $n$,
and set
\[
S = \{(a,b) \in \ZZ \times \ZZ: 0 \leq a \leq n, \quad 0 \leq b \leq n-a\},
\]
with the partial order given by
\[
(a',b') < (a,b) \quad \Longleftrightarrow \quad 
a' < a \quad \mbox{and} \quad b \leq b' \leq b + a - a'.
\]
Choose $c_0, \dots, c_{n-1} \in K$ of nonnegative valuation, 
and define a recurrence on $S$ by
\begin{align*}
f_{(0,b)} &= 1 \qquad (0 \leq b \leq n) \\
f_{(1,b)} &= c_b \qquad (0 \leq b \leq n-1) \\
f_{(a,b)} &= \frac{x_{a-1,b} x_{a-1,b+1} - 1}{x_{a-2,b+1}}
\qquad (2 \leq a \leq n, \quad 0 \leq b \leq n-a);
\end{align*}
then $g$ exists and takes values with nonnegative valuations.
Indeed, as noted in \cite{propp2}, this is basically
a special case of Dodgson condensation: 
the $f_{(a,b)}$
are connected minors of the tridiagonal
matrix
\[
\begin{pmatrix}
c_0 & 1 & 0 & & 0 & 0 \\
1 & c_1 & 1 & \cdots & 0 & 0 \\
0 & 1 & c_2 & &0 & 0 \\
& \vdots & & \ddots & & \vdots \\
0 & 0 & 0 & & c_{n-2} & 1 \\
0 & 0 & 0 & \cdots & 1 & c_{n-1} 
\end{pmatrix},
\]
and while one cannot condense this matrix (as some of 
the other connected minors vanish), one can recover the number frieze
recurrence by
instead condensing the matrix
\[
A_{ij} = \begin{cases} c_{i-1} & i = j \\
t^{(|i-j|)(|i-j|+1)/2} & i \neq j, \\
\end{cases}
\]
where $t$ is an indeterminate, then setting $t=0$ in the resulting 
polynomials.

\begin{theorem*}
The number frieze recurrence $f_{(a,b)}$ exhibits Robbins stability.
\end{theorem*}
Note that the proof will actually yield a stronger result, as in
the trivial case ($Q_s = 1$ for all $s$): 
it effectively shows that as long as the
projected precision loss is strictly less than $N$, Robbins stability
holds even using fixed point arithmetic (i.e., working
modulo $p^N$) instead of floating point arithmetic. 
\begin{proof}
Let $g'$ be an $N$-perturbation. (To simplify notation, we write 
$g(a,b)$ and $g'(a,b)$ instead of $g((a,b))$ and $g'((a,b))$.)
We prove by induction on $a$ that as long as
$r_{(a,b)}(g') < N$, we have $v(g'(a,b) - g(a,b)) \geq N - r_{(a,b)}(g')$
(and hence $v(g'(a,b)) \geq 0$, since $g(a,b)$ is known to have
nonnegative valuation); 
this gives precisely the Robbins stability bound.

Before continuing, we introduce another notational convention.
Put $r = r_{(a,b)}(g')$, and write $Y \equiv Z$ to mean
$v(Y-Z) \geq N-r$ (so in particular any star is congruent to 0). 
Note that the congruences $Y \equiv Z$ and $Y' \equiv Z'$
imply that $Y + Z \equiv Y' + Z'$ always; if $Y,Z,Y',Z'$
have nonnegative valuation, the congruences also imply that
$YY' \equiv ZZ'$. Moreover, if $Y \equiv Z$ and $Y,Z$ both
have valuation 0, then $Y^{-1} \equiv Z^{-1}$.

We now return to the induction.
For $a=0,1$, the desired inequality holds by default
because $(a,b)$ is initial.
For $a=2$, the denominator of $f_{(a,b)}$ is $x_{(0,b+1)}$,
and $g'(0,b+1) = g(0,b+1) + * = 1 + *$ has valuation 0,
so again the desired inequality follows.
For $a=3$ and $0 \leq b \leq n-3$, we have
\begin{align*}
g(3,b) &= \frac{g(2,b) g(2,b+1) - 1}{g(1,b+1)} \\
g'(3,b) &= \frac{(1 + *) g'(2,b) g'(2,b+1) - (1 + *)}{(1 + *)g'(1,b+1)};
\end{align*}
by the induction hypothesis, $g'(2,b) = g(2,b) + *$,
$g'(2,b+1) = g(2,b+1) + *$, and $g'(1,b+1) = g(1,b+1) + *$, so
\[
g'(3,b) = \frac{g(2,b) g(2,b+1) - 1 + *}{g(1,b+1) + *}.
\]
Since $Q_{(a',b')}(g') = 1$ for $a' = 0,1$, and
since for $a'=2$ we have as above $Q_{(a',b')}(g') = 1 + *$, we have
\begin{align*}
r &= \max_{(a',b') \leq (a,b)} \{ v(Q_{(a',b')}(g')) \} \\
&= v(Q_{(a,b)}(g')) \\
&= v(g'(1,b+1)).
\end{align*}
Hence (since $r<N$ by assumption) we have $g'(1,b+1) < N$, yielding
$v(g'(1,b+1) + *) = v(g'(1,b+1))$; in particular,
$v(g(1,b+1)) = v(g'(1,b+1)) = r$.
We can now write
\begin{align*}
g'(3,b) &= \frac{g(2,b) g(2,b+1) - 1 + *}{g(1,b+1) + *} \\
 &= \frac{((g(2,b) g(2,b+1) - 1)/g(1,b+1)) + (*/g(1,b+1))}{1 + */g(1,b+1)} \\
&= \frac{g(3,b) + (*/g(1,b+1))}{1 + (*/g(1,b+1))} \\
&\equiv g(3,b),
\end{align*}
as desired.

Suppose now that $a \geq 4$, $r_{(a,b)}(g') < N$, and the
induction hypothesis holds for all pairs $(a',b') < (a,b)$;
in particular, we have $v(g'(a',b')) \geq 0$ whenever $(a',b') < (a,b)$.
To eliminate some indices, put
\begin{gather*}
A = g(a-4,b+2), \\
B = g(a-3,b+1), \quad C = g(a-3,b+2), \\
D = g(a-2,b), \quad E = g(a-2,b+1), \quad F = g(a-2,b+2) \\
G = g(a-1,b), \quad H = g(a-1,b+1), \\
I = g(a,b)
\end{gather*}
and likewise with primes; note that $A, \dots, I$ all have nonnegative
valuation, as do $A', \dots, H'$ by the induction hypothesis.
We then have
\begin{gather*}
E' = \frac{B'C' - 1 + *}{A' + *}, \\
G' = \frac{D'E' - 1 + *}{B' + *}, \qquad H' = \frac{E'F' - 1 + *}{C' + *}, \\
I' = \frac{G'H' - 1 + *}{E'+ *},
\end{gather*}
because $g'$ is an $N$-perturbation and
$v(g'(a',b')) \geq 0$ for $a' < a$. (More explicitly, the definition
of an $N$-perturbation implies that $E' = (B'C'(1+*) - (1+*))/(A'(1+*))$
and the like, but the product of each lettered quantity with a star
is again a star.)
We also have four analogous equations without the primes and stars.
Moreover, if $(a',b') < (a,b)$, we have $r \geq r_{(a',b')}(g')$ by the
way the projected precision loss is defined, so the induction
hypothesis implies in particular that $g'(a',b') \equiv g(a,b)$;
in particular, we have
\[
A' \equiv A, \dots, H' \equiv H,
\]
and we wish to show that $I' \equiv I$.

By the induction hypothesis, we have $v(E') \geq 0$.
If $v(E') = 0$, then $G' \equiv G, H' \equiv H, E' \equiv E$ imply
$G'H' - 1 + * \equiv GH - 1$ and $E'+* \equiv E$. Since
$N > r$, the congruence $E' \equiv 
E$ and the assumption $v(E') = 0$ imply $v(E) = 0$, 
and so $(E' + *)^{-1} \equiv E^{-1}$. 
Consequently
\begin{align*}
I' &= \frac{G'H'-1+*}{E'+*} \\
&\equiv \frac{GH-1}{E} \\
&= I
\end{align*}
as desired.

Since the case $v(E) =0$ is okay, we assume hereafter that $v(E') > 0$; then
$v(B'C' - 1 + *) > 0$, and hence $v(B'C' - 1) > 0$. 
Since $v(B') \geq 0$, $v(C') \geq 0$, and
$0 = v(1) \geq \min\{v(B'C'), v(1-B'C')\}$,
this is only possible if
$v(B') = v(C') = 0$.

We now compute
\begin{align*}
I' &= \frac{G'H' - 1 + *}{E' + *} \\
&= \frac{(D'E'-1 + *)(E'F'-1 + *) - (B'+*)(C'+*)(1+*)}{(B'+*)(C'+*)(E'+*)} \\
&= \frac{D'E'E'F' - D'E' - E'F' + 1 - B'C' + *}{B'C'E' + *} \\
&= \frac{D'E'E'F' - D'E' - E'F' - A'E' + *}{B'C'E' + *} \\
&= \frac{D'E'F' - D' - F' - A' + (*/E')}{B'C' + (*/E')}.
\end{align*}
As before, we have $D'E'F' \equiv DEF$, $D' \equiv D$, $F' \equiv F$,
$A'\equiv A$, and $B'C' \equiv BC$.
Moreover, from the definition of the projected precision loss, we have
\begin{align*}
r &= \max_{(a',b') \leq (a,b)} \{ v(Q_{(a',b')}(g')) \} \\
&\geq v(Q_{(a,b)}(g')) \\
&= v(E'),
\end{align*}
and so $*/E' \equiv 0$.

Since $r < N$, the facts that $v(B'C') = 0$ and
$B'C' \equiv BC$ together imply that $v(BC) = 0$;
then the congruence
$BC \equiv B'C' + (*/E')$ implies
$(B'C' + (*/E'))^{-1} \equiv (BC)^{-1}$.
This together with the previous mentioned congruences and the equation
\[
I = \frac{DEF - D - F - A}{BC}
\]
yields
$I' \equiv I$, as desired.
\end{proof}

Note that in this example,
the precision bound given by Robbins stability is not always
sharp if one fixes $(a,b)$ and varies over all $N$-perturbations. 
For instance,
for $K = \QQ$ with the $3$-adic valuation, take
\[
(c_0, \dots, c_5) = (1, 3^m - 1, -1, 1, -11, 22).
\]
For $m$ and $N$ sufficiently large (say $m > 5$ and $N \geq 2m$),
the projected precision loss is $m$ (achieved by $g(1,1) = -3^m$),
but experiments suggest that $v(g'(5,0)-g(5,0)) \geq N - m + 5$ always.
It would be interesting to find a more precise version of the
projected precision loss that detects such ``localized
disruptions'', specifically by
relaxing the restriction that the bound can only get worse with
each successive term. Such a formulation of the stability phenomenon
may even suggest progress towards Robbins's original conjecture
or generalizations.

Although all our examples have been recurrences over $\QQ$, with $v$
equal to a $p$-adic valuation, we have taken care to make our setup
more general. In particular, one could use our framework to look
at Robbins stability in $\QQ(x)$, with $v$ the $x$-adic valuation.
This might serve as a bridge between the Laurent phenomenon and Robbins
stability.

We conclude by mentioning some further experiments the first author
has conducted with Punyashloka Biswal. Namely, we have been applying Robbins's
testing regimen to other recurrences exhibiting the
``Laurent phenomenon'' of \cite{fz}: compute pairs of
$N$-perturbations using $N$-digit $p$-adic floating point arithmetic
(generating the undetermined $p$-adic digits at random),
and compare their difference to the projected precision loss
predicted by Robbins stability. (This is somewhat easier than comparing one
$N$-perturbation to the exact solution.)
Two families of examples we have considered, which both appear to exhibit
Robbins stability, are the Somos sequences
\[
x_0 = x_1 = \cdots = x_{k-1} = 1, \qquad
x_{n+k} = 
\frac{\sum_{1 \leq i \leq \lfloor k/2\rfloor} a_i x_{n+i} x_{n+k-i}}{x_n}
\]
for $k=4, 5, 6, 7$, and the sequences
\[
x_{n+2} = \frac{x_{n+1}^2 + c x_{n+1} + d}{x_n}
\]
given in \cite[Example~5.4]{fz}.
Notably, the latter example seems to require the correction
term $\min\{0, v(g(s))\}$ that we introduced into the definition of
Robbins stability.

\subsection*{Acknowledgments}

Thanks to Joe Buhler for referring us to the formulation of Robbins's 
conjecture appearing here.
The first author is supported by NSF grant DMS-0400727,
and the second author is supported by NSA grant H92830-04-1-0054.

\end{document}